\documentclass[12pt, leqno]{article}
\usepackage{amssymb}
\def\q \m#1#2{{\raise1pt\hbox{$#1$}\kern-1pt\big/
               \kern-1pt\raise-1pt\hbox{$#2$}}}

\def\bZ{{\rm \bf Z}}
\def\bQ{{\rm \bf Q}}
\def\bC{{\rm {\bf C}}}

\def\bH{{\rm \bf H}}
\def\sR{{ \rm \scriptsize  \bf R}}

\def\sC{{ \rm \scriptsize  \bf C}}

\newfam\msbfam
\font\twelmsb=msbm10 at 12pt
\font\tenmsb=msbm10 at 10 pt
\font\sevenmsb=msbm10 at 7pt
\font\fivemsb=msbm10 at 5pt

\textfont\msbfam=\twelmsb
\textfont\msbfam=\tenmsb
\scriptfont\msbfam=\sevenmsb
\scriptscriptfont\msbfam=\fivemsb

\begin{document}

\renewcommand{\theequation}{\thesection.\arabic{equation}}
\setcounter{equation}{0}

\centerline{\Large {\bf On Elliptic Genera and Foliations}}

\vskip 8mm  
\centerline{\bf Kefeng Liu,\footnote{Partially supported by the Sloan 
Fellowship and an NSF grant.}  Xiaonan Ma\footnote{Partially supported by SFB 288.}
  and  Weiping Zhang\footnote{Partially supported by NSFC, MOEC and the Qiu Shi Foundation.}}
\vskip 5mm

{\small {\bf Abstract.} We prove several vanishing 
theorems for a class of generalized elliptic genera 
on foliated manifolds, by using
classical equivariant index theory. The 
main techniques are the use of the Jacobi theta-functions and the 
construction of  a new class of 
elliptic operators associated to  foliations.}

\vspace{.15in}

{\bf \S 1. Introduction}

\vspace{.15in}

The main purpose of this paper is to prove some vanishing theorems 
 of characteristic numbers for foliated manifolds with group actions.
 Such type of results are usually proved by using 
index theorems or fixed point theorems 
for foliations as developed by 
Connes, Skandalis, Heitsch, Lazarov and so on (cf. [Co] and [HL1]). In this paper we 
take a rather different route. Instead of using the heavy machinery about index theory for 
foliations as developed by various people, we 
use certain new  elliptic operators particularly designed to 
study integrable subbundles with spin struture, the so-called sub-Dirac 
operator. With the help of such operators, we are able to prove our theorems
for foliated manifolds by using the classical index theory (Compare with [HL2]). 
Such operators are also  used in [LiuZ] to prove certain 
vanishing theorem for characteristic numbers of foliations
with spin leaves of positive scalar curvature, which were previously proved 
by Connes (cf. [Co, pp. 182 and 270]) by using noncommutative index theory.  

More precisely,  let $(M,F)$ be 
a transversally oriented compact foliated manifold such 
that the integrable bundle $F$ is spin
and carries a fixed spin structure.
Assume that there is an effective $S^3$-action on $M$ which preserves the 
leaves induced by $F$ and also the spin structure on $F$. Then 
a special case of our result shows that the Witten genus 
([W, (17)]) of $M$ vanishes if the first Pontryagin class of $F$ verifies $p_1(F)=0$. 
Note that here we do not assume that the manifold $M$ is spin, so that the Witten genus
under consideration in general might not be an integer.

On the other hand, we will also prove certain vanishing theorems for spin manifolds 
with split tangent bundles, by using the similar technique. The elliptic 
genera we derive in this situation can be viewed as 
interpolations between the various classical 
elliptic genera. They are actually the mixture of the two universal 
elliptic genera. It is interesting to note that, under some mild 
conditions,  we get quite general vanishing theorems. Similar theorems can be 
proved for loop group representations. 

This paper is organized  as follows. In Section 2, we introduce 
the sub-Dirac operator. In Section 3, we state our 
main vanishing theorem for elliptic genus on foliations, which will be 
proved in Section 4 by combining the construction in 
Section 2 with Jacobi-theta functions.
In Section 5, we prove several vanishing theorems for certain twisted elliptic
genera associated to spin manifolds with split tangent bundle.
In Section 6, we point out other generalizations and 
state a conjecture concerning the 
vanishing of the Witten genus of a foliation with spin leaves of
positive Ricci curvature. 

\vspace{.15in}

{\bf \S 2. Sub-Dirac operator}

\vspace{.15in}

Let $M$ be an even dimensional smooth compact oriented manifold. 
Let $F$ be a sub-bundle of the 
tangent vector bundle $TM$ of $M$. Let $g^{TM}$ be a Riemannian metric on $TM$.
Let $F^\bot$ be the orthogonal complement to $F$ in $TM$. Then one has the 
orthogonal splittings
$$TM=F\oplus F^\bot,$$
$$g^{TM}=g^F\oplus g^{F^\bot}.\eqno (2.1)$$
Moreover, one has the obvious identification that 
$$TM/F\simeq F^\bot.\eqno(2.2)$$

{}From now on we make the special assumption that $F$ is  even dimensional,
oriented, spin and carries 
a fixed spin structure. Then $F^\bot$ carries an induced orientation. 
Set $2p=\dim F$ and $2r=\dim F^\bot$.

Let $S(F)$ be the bundle of spinors associated to $(F, g^F)$. For any 
$X\in F$, denote by $c(X)$ the Clifford action of $X$ on $S(F)$. 
We have the splitting 
$$S(F)=S_+(F)\oplus S_-(F)\eqno(2.3)$$
and $c(X)$ exchanges $S_\pm(F)$.

Let $\Lambda(F^{\bot,*})$ be the exterior algebra bundle of $F^\bot$. 
Then $\Lambda(F^{\bot,*})$ carries a canonically induced metric 
$g^{\Lambda(F^{\bot,*})}$ from $g^{F^\bot}$.
By using $g^{F^\bot}$, one has the canonical identification 
$F^\bot\simeq F^{\bot, *}$. For any $U\in F^\bot$, 
let $U^*\in F^{\bot, *}$ be the corresponding dual of 
$U$ with respect to $g^{F^\bot}$.
Now for $U\in F^\bot$, set 
$$ c(U)=U^*\land -i_U, \ \widehat{c}(U)=U^*\land+i_U,\eqno(2.4)$$ 
where $U^*\land $ and $i_U$ are the exterior and inner multiplications by 
$U^*$ and $U$ on $\Lambda(F^{\bot,*})$ respectively. One has the following 
obvious identities,
$$c(U)c(V)+c(V)c(U)=-2\langle U,V\rangle_{g^{F^\bot}},$$ 
 $$\widehat{c}(U)\widehat{c}(V)+\widehat{c}(V)\widehat{c}(U)
=2\langle U,V\rangle_{g^{F^\bot}},$$
$$c(U)\widehat{c}(V)+\widehat{c}(V)c(U)=0  \eqno(2.5)$$
for $  U,\ V\in F^\bot$.

Let $h_1,\cdots, h_{2r}$ be an oriented local orthonormal basis of $F^\bot$. 
Set 
$$\tau\left(F^\bot, g^{F^\bot}\right)=(\sqrt{-1})^{r}c(h_1)\cdots c(h_{2r}).
\eqno(2.6)$$
Then clearly 
$$\tau\left(F^\bot, g^{F^\bot}\right)^2={\rm Id}_{\Lambda(F^{\bot,*})}.
\eqno(2.7)$$

Denote 
$$\Lambda_\pm\left(F^{\bot,*}\right)=\left\{ h\in \Lambda\left(F^{\bot,*}\right):\ 
\tau\left(F^\bot, 
g^{F^\bot}\right)h=\pm h \right\}.\eqno(2.8)$$
Then $\Lambda_\pm(F^{\bot,*})$ are sub-bundles of $\Lambda(F^{\bot,*})$. 
Also, one verifies that for any $h\in F^\bot$, 
$c(h)$ anticommutes with $\tau$, while $\widehat{c}(h)$ commutes with $\tau$. 
Thus $c(h)$ exchanges $\Lambda_\pm(F^{\bot,*})$.

We will view both vector bundles 
$$S(F)=S_+(F)\oplus S_-(F)\eqno(2.9)$$ 
and 
$$\Lambda\left(F^{\bot,*}\right)=\Lambda_+\left(F^{\bot,*}\right)\oplus
 \Lambda_-
\left(F^{\bot,*}\right)\eqno(2.10)$$ 
as super-vector bundles. 
Their ${\bf Z}_2$-graded tensor product is given by 
$$S(F)\widehat{\otimes}\Lambda\left(F^{\bot,*}\right)=\left[S_+(F)\otimes\Lambda_+
\left(F^{\bot,*}\right)
\oplus S_-(F)\otimes\Lambda_-\left(F^{\bot,*}\right)\right]$$
$$\bigoplus\left[S_+(F)\otimes\Lambda_-\left(F^{\bot,*}\right)
\oplus S_-(F)\otimes\Lambda_+\left(F^{\bot,*}\right)\right].\eqno(2.11)$$
For $X\in F$, $U\in F^\bot$, the operators $c(X)$, 
$c(U)$ and $\widehat{c}(U)$ extends naturally to 
$S(F)\widehat{\otimes}\Lambda(F^{\bot,*})$.

The connections $\nabla^F$, $\nabla^{F^\bot}$ lift to $S(F)$ and 
$ \Lambda(F^{\bot,*})$ naturally, and preserve the splittings 
(2.9) and (2.10). We write them as
$$\nabla^{S(F)}=\nabla^{S_+(F)}\oplus\nabla^{S_-(F)},\ \ \ 
\nabla^{\Lambda(F^{\bot,*})}=\nabla^{\Lambda_+(F^{\bot,*})}
\oplus\nabla^{\Lambda_-(F^{\bot,*})}.\eqno(2.12)$$
Then $S(F)\widehat{\otimes}\Lambda(F^{\bot,*})$ carries the induced 
tensor product connection
$$\nabla^{S(F)\widehat{\otimes}\Lambda(F^{\bot,*})}=\nabla^{S(F)}\otimes 
{\rm Id}_{\Lambda(F^{\bot,*})}+{\rm Id}_{S(F)}\otimes\nabla^{\Lambda(F^{\bot,*})} .
\eqno(2.13)$$
And similarly for $S_\pm(F)\widehat{\otimes}\Lambda_\pm(F^{\bot,*})$.

Let $S\in \Gamma(T^*M \otimes \mbox{End}(TM))$ be defined by 
$$\nabla^{TM}=\nabla^{F}+\nabla^{F^\bot}+S.\eqno(2.14)$$ 
Then for any $X\in TM$, 
$S(X)$ exchanges $F$ and $F^\bot$ and is skew-adjoint with
 respect to $g^{TM}$.

For any vector bundle $E$ over $M$, by an integral polynomial of $E$ we 
will mean a vector bundle $\varphi(E)$ which is a polynomial in the exterior and
 symmetric powers of $E$ with integral coefficients.

Let $\psi(F)$ (resp. $\varphi(F^\bot)$)
 be an integral polynomial of $F$ (resp. $F^\bot$), then $\psi(F)$ 
(resp. $\varphi(F^\bot)$)
carries a naturally induced metric $g^{\psi(F)}$ 
(resp. $g^{\varphi(F^\bot)}$) from $g^F$ (resp. $g^{F^\bot}$) and also 
a naturally induced Hermitian connection $\nabla^{\psi(F)}$ 
(resp. $\nabla^{\varphi(F^\bot)}$) induced from 
$\nabla^F$ (resp. $\nabla^{F^\bot}$).

Our main concern will be on the $\bZ_2$-graded vector bundle 
$$\left(S(F)\widehat{\otimes}\Lambda\left(F^{\bot,*}\right)\right)\otimes
\psi(F)\otimes \varphi\left(F^\bot\right)\eqno(2.15)$$ 
which is 
$$\left[S_+(F)\otimes\Lambda_+\left(F^{\bot,*}\right)\otimes\psi(F)\otimes 
\varphi\left(F^\bot\right) 
\oplus S_-(F)\otimes\Lambda_-\left(F^{\bot,*}\right)\otimes\psi(F)\otimes \varphi
\left(F^\bot\right)\right]$$
$$\bigoplus\left[S_+(F)\otimes\Lambda_-\left(F^{\bot,*}\right)\otimes
\psi(F)\otimes \varphi\left(F^\bot\right) 
\oplus S_+(F)\otimes\Lambda_-\left(F^{\bot,*}\right)
\otimes\psi(F)\otimes \varphi\left(F^\bot\right)\right].$$

The Clifford actions $c(X)$, $c(U)$ and $\widehat{c}(U)$ for $X\in F$, 
$U\in F^\bot$ extend further to these bundles by acting as identity
 on $\psi(F)\otimes \varphi(F^\bot)$.

We can also form the tensor product metric on the new bundles as well as 
the tensor product connection on
$(S(F)\widehat{\otimes}\Lambda(F^{\bot,*}))\otimes \psi(F)\otimes \varphi(F^\bot) $
given by
$$\nabla^{(S(F)\widehat{\otimes}\Lambda(F^{\bot,*}))\otimes \psi(F)\otimes \varphi(F^\bot) }
=\nabla^{S(F)\widehat{\otimes}\Lambda(F^{\bot,*})}\otimes 
{\rm Id}_{\psi(F)\otimes \varphi(F^\bot)}+
{\rm Id}_{S(F)\widehat{\otimes}\Lambda(F^{\bot,*})}
\otimes \nabla^{\psi(F)\otimes \varphi(F^\bot)},\eqno(2.16)$$
where $\nabla^{\psi(F)\otimes \varphi(F^\bot)}$ is the tensor product connection on
$\psi(F)\otimes \varphi(F^\bot)$ obtained from $\nabla^{\psi(F)}$ and
$\nabla^{\varphi(F^\bot)}$,
as well as on the $\pm$ subbundles.

Now let $\{ f_i\}^{2p}_{i=1}$ be an oriented orthonormal 
basis of $F$. Recall 
that $\{ h_s\}^{2r}_{s=1}$ is an oriented orthonormal basis of $F^\bot$. 
The elliptic 
operators which are the main concern of this section can be defined as 
follows. It is introduced mainly for the reason that the vector bundle
$F^\bot$ might well be non-spin.

$\ $

{\bf Definition 2.1.} Let $D_{F, \psi(F)\otimes \varphi(F^\bot)}$ be the 
operator which maps
$\Gamma(S(F)\widehat{\otimes}\Lambda(F^{\bot,*})\otimes \psi(F)\otimes \varphi(F^\bot))$ 
to itself defined by  
$$D_{F, \psi(F)\otimes \varphi(F^\bot)}= \sum^{2p}_{i=1}c(f_i)\nabla_{f_i}^{(S(F)
\widehat{\otimes}\Lambda(F^{\bot,*}))\otimes \psi(F)\otimes \varphi(F^\bot)}
 +\sum^{2r}_{s=1}c(h_s)
\nabla_{h_s }^{(S(F)\widehat{\otimes}\Lambda(F^{\bot,*}))\otimes \psi(F)\otimes \varphi(F^\bot)}$$
$$+{1\over 2}\sum_{i,j=1}^{2p}\sum^{2r}_{s=1}\langle
 S(f_i)f_j, h_s\rangle c(f_i)c(f_j)c(h_s)
 +{1\over 2}\sum_{s,t=1}^{2r}\sum^{2p}_{i=1}
\langle S(h_s)h_t,f_i\rangle c(h_s)c(h_t)c(f_i).\eqno(2.17)$$

$\ $

It is easy to verify that $D_{F, \psi(F)\otimes \varphi(F^\bot)}$ is a first order 
formally self adjoint elliptic differential operator. Furthermore, it 
anticommutes with the ${\bf Z}_2$ grading operator of the super vector bundle 
$(S(F)\widehat{\otimes}\Lambda(F^{\bot,*}))\otimes \psi(F)\otimes \varphi(F^\bot)$.

Let $D_{F, \psi(F)\otimes \varphi(F^\bot), +}$ (resp. 
$D_{F,\psi(F)\otimes  \varphi(F^\bot), -}$) be the 
restriction of $D_{F, \psi(F)\otimes \varphi(F^\bot)}$ to the even (resp. odd) subbundle 
of $(S(F)\widehat{\otimes}\Lambda(F^{\bot,*}))\otimes \psi(F)\otimes \varphi(F^\bot)$. 
Then one has 
$$D_{F, \psi(F)\otimes \varphi(F^\bot), +}^*=
D_{F, \psi(F)\otimes \varphi(F^\bot), -}.\eqno(2.18)$$

\vspace{.15in}

{\bf Remark 2.2.}  Locally, $D_{F,\psi(F)\otimes \varphi(F^\bot)}$ may be seen as a 
twisted Dirac operator. The key point {\em here} is that its definition 
relies {\em only} on the spin structure of $F$
(Compare with Remark 4.4 in Section 4).

$\ $

The following result follows easily from the Atiyah-Singer index theorem [AS].

\vspace{.15in}

{\bf Theorem 2.3.} {\em The following index formula holds,}
$${\mbox{ind}} \left( D_{F, \psi(F)\otimes \varphi(F^\bot), +} \right)
=\left\langle\widehat{A}(F){ \mbox{ch}}
(\psi(F)) L\left(F^\bot\right)
{ \mbox{ch}}\left(\varphi\left(F^\bot\right)\right) , [M]
\right\rangle.\eqno(2.19)$$  

\vspace{.15in}

Now assume that $M$ admits an $S^1$-action which preserves 
$g^{TM}$, as well as the spin structure on $F$. Then 
it also preserves the splittings in (2.1). 
Furthermore, an equivariant version of the  index formula (2.19) still holds. 

More precisely, let $\{N\}$ be the set of connected components of the fixed 
point set of this circle action. Assume 
that when restricted to the fixed point set, we have the equivariant 
decompositions
$$F|_N = F_0\oplus \left(\oplus_j E_j\right),\ \ \ \
F^\bot =F_0^\bot \oplus \left(\oplus_j L_j\right)\eqno(2.20)$$
such that the generator $e^{2\pi it}\in S^1$ acts trivially on 
the real vector bundles $F_0$ and $F_0^\bot$, 
and acts on the complex vector bundles $E_j$ and $L_j$ by multiplications by
 $e^{2\pi im_j t}$ and $e^{2\pi in_j t}$  respectively. 
Let $\{2 \pi i x_j^k \} $ be the  Chern roots of $E_j$ and 
$\{2 \pi i z_j^k \}$ be the  Chern roots of $L_j$. 
(Note that in our notation, if $E_j$ is a complex line bundle and $R^{E_j}$
is the curvature of a connection on $E_j$, then ${-R^{E_j}\over 2 \pi i }= 2 \pi i x_j$.)

By (2.20), $F_0,\ F^{\bot}$ are naturally oreinted. We fixe the orientation 
on $N$ induced by the orientations on  $F_0,\ F^{\bot}$.

The following result follows easily from the equivariant index theorem of 
Atiyah, Bott, Segal and Singer (cf. [AS]).

\vspace{.15in}

{\bf Theorem 2.4} {\em  The following  equivariant 
index formula for the Lefschetz number of the generator $g=e^{2\pi it}\in S^1$ 
associated to the elliptic operator 
$D_{F, \psi(F)\otimes \varphi(F^\bot), +}$ holds,
$$L(g)=
\sum_N \left\langle \widehat{A}(F_0) L\left(F_0^\bot\right) A(F, t)L\left(F^\bot, t\right)
{\mbox{\rm ch}}_g\left(\psi(F|_N)\right){\mbox{\rm ch}}_g\left(\varphi\left(F^\bot|_N
\right)\right), [N]
\right\rangle,\eqno(2.21)$$
where 
$$A(F, t) = \prod_{j,k}\frac{1}{2{ \mbox{\rm sinh}}(\pi i (x_j^k +m_jt))},\ \ \ 
L\left(F^\bot, t\right)
=\prod_{j,k} \frac{1}{{ \mbox{\rm tanh}}(\pi i (x_j^k+n_jt))} \eqno(2.22)$$
and ${\mbox{\rm ch}}_g$ denotes the 
equivariant Chern character, for examples,
 $${\mbox{\rm ch}}_g(F|_N) ={\mbox{\rm ch}}(F_0)+ \sum_{j,k} e^{2 \pi i (x_j^k+m_jt)}, \ \ \ 
{\mbox{\rm ch}}_g\left(F^\bot|_N \right)={\mbox{\rm  ch}}
\left(F^\bot_0\right)+ \sum_{j,k} e^{2 \pi i (z_j^k+n_jt) }.
\eqno(2.23)$$}

{\bf \S 3. Elliptic genus for foliations}

$\ $

For any vector bundle $E$, let us denote the two operations in $K$-theory, 
the total symmetric and exterior power operations,  by
$${\rm Sym}_q(E)=1+qE+q^2{\rm Sym}^2( E)+\cdots, $$
$$\Lambda_q(E)=1+tE+t^2\Lambda^2(E)+\cdots,\eqno(3.1)$$
where $q$ is a parameter. Recall the relations:
$${\rm Sym}_q(E_1-E_2)={\rm Sym}_q(E_1)\cdot \Lambda_{-q}(E_2); 
\ \ \ \Lambda_q(E_1-E_2)=\Lambda_q(E_1)\cdot {\rm Sym}_{-q}(E_2).$$

In what follows, we will take $\psi(F)$ to be
the Witten element [W]
$$\Psi_q(F)=\otimes_{j=1}^\infty {\rm Sym}_{q^j}(F-{\mbox{dim}}\, F)\eqno(3.2)$$
with $q$ a parameter. 

We now further assume in this section that $F$ is a nontrivial integrable sub-bundle
of $TM$. Then $F$ induces a foliation on $M$. We make the basic assumption 
that 
the $S^1$-action on $M$ preserves the leaves induced by $F$.

Recall that the equivariant cohomology group $H^*_{S^1}(M)$ is defined to be 
the usual cohomology group of the space $ES^1\times_{S^1} M$, 
where $ES^1$ denotes 
the universal principal $S^1$-bundle over the classifying space $BS^1$. Here 
we take cohomology with rational coefficients. The projection 
$$\pi: \ ES^1\times_{S^1} M\rightarrow B{S^1}\eqno(3.3)$$ 
induces a map 
$$\pi^*: \ H^*_{S^1}({\rm pt.})\rightarrow H^*_{S^1}(M)\eqno(3.4)$$ 
which  makes 
$H^*_{S^1}(M)$ a module over $H^*_{S^1}({\rm pt.})\simeq {\bf Q}[[u]]$ 
with $u$ a generator of degree $2$.

Let $p_1(F)_{S^1}$ be the equivariant first Pontryagin class of $F$. We can now state
the main result of this section as follows.

\vspace{.15in}

{\bf Theorem 3.1.} {\em If $M$ is connected and the $S^1$-action on $M$ 
is nontrivial and $p_1(F)_{S^1}=n\cdot \pi^* u^2$ for some 
integer $n$, then the equivariant index of 
$D_{F, \Psi_q(F)\otimes \varphi(F^\bot), +}$ 
is $0$. As a consequence, for any Pontryagin class $p(TM/F)$ of $TM/F$, 
we have 
$$\left\langle\widehat{A}(F){\mbox{\rm ch}}\left(\Psi_q(F)\right) p(TM/F), [M]
\right\rangle=0.\eqno(3.5)$$
In particular, the Witten genus [W] of $M$, which is defined by
$$\left\langle \widehat{A}(TM){\mbox{\rm ch}}\left(\Psi_q(TM)\right), [M]
\right\rangle, \eqno(3.6)$$
vanishes. }

\vspace{.15in}

If the $S^1$-action is induced from an effective $S^3$-action which also 
preserves the foliation and the spin structure on $F$, then 
one can show that $p_1(F)_{S^1}=n\cdot \pi^* u^2$
is equivalent to the condition that $p_1(F)=0$. This gives 
us the following

\vspace{.15in}

{\bf Corollary 3.2.} {\em Assume that $M$ is connected and that 
there is an effective $S^3$-action that preserves 
the foliation and the spin structure on $F$, and that $p_1(F)=0$, 
then  the equivariant index of 
$D_{F,\Psi_q(F)\otimes \varphi(F^\bot), +}$ is $0$. In particular,
the vanishing formula (3.5) holds and the Witten genus 
given by (3.6) vanishes.}

\vspace{.15in}

Theorem 3.1 and Corollary 3.2 will be proved in the next section.

\vspace{.15in}

{\bf \S 4. Proof of Theorem 3.1}

\vspace{.15in}

Let us first recall the defintion of the Jacobi-theta functions 
as in [Ch].

For $v\in \bC$, $\tau \in \bH= \{ \tau \in \bC, {\rm Im} \tau >0\}$,
 $q=e^{2\pi i\tau}$, let 
 $\theta(v,\tau)$ denote the classical Jacobi theta-function
$$\theta(v, \tau) =c(q)q^{1/8}2{\mbox{sin}}\, (\pi v) \prod_{n=1}^\infty
\left(1-q^n e^{2\pi iv}\right)\left(1-q^n e^{-2\pi iv}\right)\eqno(4.1)$$ 
where $c(q) = \prod_{n=1}^\infty\left(1-q^n\right)$.  Set
$$\theta'(0, \tau) =\left. 
\frac{\partial \theta(v, \tau) }{\partial v}\right|_{v=0}.
\eqno(4.2)$$

Since any Pontryagin class $p(F^\bot)$ of $F^\bot$
can be written as a  linear 
combination with rational coefficients of  classes of the form 
$L(F^\bot){\rm ch}(\varphi(F^\bot))$, we can and we will assume first that
$p(F^\bot)$ is of homogeneous degree $2l$ and that $\varphi(F^\bot)$ verifies that
$L(F^\bot){\rm ch}(\varphi(F^\bot))$ is a nonzero rational multiple of $p(F^\bot)$.

Let $g=e^{2\pi it}\in S^1$ be a generator of the $S^1$-action.
Let $\{N\}$ denote the set of connected components of the fixed point set
of the $S^1$-action. Since the $S^1$-action preserves the leaves induced by $F$, according to 
Lemma 3.3 in [HL2], it induces a trivial action on $F^\bot|_N$.
Assume that the bundle $F|_N$ has the decomposition 
$$F|_N = F_0\oplus \left(\oplus_j E_j\right),\eqno(4.3)$$ 
where each $E_j$ is a complex vector bundle on which 
$e^{2\pi it}$ acts by $e^{2\pi i m_j t}$,
while the $S^1$ acts trivially on the real vector bundle $F_0$.

Let $\{ 2 \pi i x_j^k \} $ denote the  Chern roots of $E_j$, 
and let $\{ \pm 2 \pi i  y_j\}$ denote the Chern roots of $F_0\otimes_\sR \bC$.
By Theorem 2.4 one deduces easily that the Lefschetz number $L(g)$ associated to the
operator $D_{F, \Psi_q(F)\otimes \varphi(F^\bot), +}$ is given by
$$H(t, \tau)= (2 \pi i)^{-p} \sum_N \left\langle 
H(F_0, \tau)
\prod_{j,k} \left( \frac{\theta'(0, \tau)}{\theta(x_j^k +m_jt, \tau)} \right)
L\left(F^\bot|_N\right) 
{\mbox{\rm ch}}\left(\varphi\left(F^\bot|_N\right)\right), [N]
\right\rangle,\eqno(4.4)$$
where the term $ H(F_0, \tau)$ denotes the characteristic class 
$$H(F_0, \tau)= \prod_j 
\left(2 \pi i y_j\frac{\theta'(0, \tau)}{\theta(y_j, \tau)}\right).\eqno(4.5)$$

Considered as function of $(t, \tau)$, we can obviously extend $H(t,\tau)$ to 
meromorphic function on $\bC \times \bH$. Note that this function
is holomorphic in $\tau$.

Recall that $2p=\dim F$ and that $L(F^\bot){\mbox{ch}}(\varphi(F^\bot))$
is of homogeneous degree $2l$. As the $S^1$ action on $M$ induces a trivial 
action on $F^{\bot}|_N$, we know that $\dim N= 2r + \dim F_0$.

\vspace{.15in}

{\bf Lemma 4.1.} {\em  The following 
formulas hold under the modular transformations,} 
$$H\left(\frac{t}{\tau}, -\frac{1}{\tau}\right)=\tau^{p+r-l} e^{-\pi i nt^2/\tau}H(t, \tau), 
\ \ \ \ 
H(t, \tau+1) =H(t, \tau).\eqno(4.6)$$

{\em Proof.} One deduces easily that the condition on $p_1(F)_{S^1}$ implies that
$$\sum_{j,k} \left(2 \pi i x_j^k +m_ju\right)^2+\sum_j (2 \pi i y_j
)^2 = n\cdot u^2,$$ 
which in turn implies that 
$$\sum_{j,k} \left(x_j^k\right)^2 +\sum_j y_j^2=0, \ \ \ 
\sum_{j,k} m_j\, x_j^k =0,\ \ \ \sum_j 
\left(\dim_\sC E_j\right) m_j^2 =n.
\eqno(4.7)$$ 

By [Ch], we have the following transformation formulas
$$
\theta \left({t \over \tau}, - {1 \over \tau}\right)= {1 \over i} \sqrt{\tau \over i}
 e^{\pi i t^2 \over \tau} \theta (t,\tau), \ \ \ \  
 \theta (t, \tau+1) = e^{ \pi i \over 4} \theta (t, \tau). \eqno(4.8)
$$

By (4.7) and (4.8), we get 
$$H\left(\frac{t}{\tau}, -\frac{1}{\tau}\right)=(2 \pi i )^{-p}\sum_N 
\left\langle H\left(F_0, -\frac{1}{\tau}\right)
\prod_{j,k} \left(\frac{\theta'\left(0, -\frac{1}{\tau}\right)}{\theta\left(x_j^k
+m_j\frac{t}{\tau},-\frac{1}{\tau}\right)}\right) L\left(F^\bot\right) 
{\mbox{ch}}\left(\varphi\left(F^\bot\right)\right), [N]\right\rangle$$
$$= (2 \pi i )^{-p} \tau^p e^{-\pi i n t^2/\tau}\sum_N  
\left\langle \prod_j 
\left(2 \pi i y_j\frac{\theta'(0, \tau)}{\theta(\tau y_j, \tau)}\right)
\prod_{j,k} \left(\frac{\theta'(0, {\tau})}{\theta
(\tau x_j^k+m_j t,{\tau})}\right)L\left(F^\bot\right) {\mbox{ch}}\left(
\varphi\left(F^\bot\right)\right), [N]\right\rangle. \eqno(4.9) $$

By comparing the $({1\over 2}\dim F_0 + r)= {1\over 2}\dim N $ homogeneous 
terms of the polynomials in 
$x$'s and $y$'s and the Chern roots of $F^{\bot}$,
 on both sides, we get the following equation: 
$$ \tau ^{-\dim F_0 /2}\left\langle \prod_j 
\left(2 \pi i  \tau y_j\frac{\theta'(0, \tau)}{\theta(\tau y_j, \tau)}\right)
\prod_{j,k}\left( \frac{\theta'(0, {\tau})}{\theta
(\tau x_j^k+m_j t,{\tau})}\right) L\left(F^\bot\right) {\mbox{ch}}\left(
\varphi\left(F^\bot\right)\right) ,[N]\right\rangle$$
$$ = \tau^{r-l} \left\langle 
H(F_0, \tau)
\prod_{j,k} \left( \frac{\theta'(0, \tau)}{\theta(x_j^k +m_jt, \tau)} \right)
L\left(F^\bot|_N\right) 
{\mbox{\rm ch}}\left(\varphi\left(F^\bot|_N\right)\right), 
[N]\right\rangle. \eqno(4.10) $$

By (4.9), (4.10), we get the first identity of (4.6). 
 By using (4.8), the second identity  can also be verified easily.  Q.E.D.

\vspace{.15in}

{\bf Lemma 4.2.} {\em  For $a,\  b \in 2{\bf Z}$, the following identity holds,}
$$H(t+a\tau +b, \tau) = e^{\pi in (a^2\tau +2a t)} H(t, \tau).\eqno(4.11)$$

{\em Proof.} By [Ch], for $a,b\in 2 \bZ$, we have the transformation formula 
for the theta-function,
$$\theta(t+a\tau +b, \tau) = e^{\pi i (a^2\tau +2a t)} \theta(t, \tau).
\eqno(4.12)$$

By using (4.7) and (4.12), 
we  obtain immediately the wanted identity.  Q.E.D.

\vspace{.15in}

Let $\bf{R}$ denote the real number field.

\vspace{.15in}

{\bf Lemma 4.3.} {\em  The function $H(t, \tau)$ is  holomorphic  
 for $(t, \tau)\in {\bf R}\times {\bf H}$.}

\vspace{.15in}

{\bf Remark 4.4.} Lemma 4.3
is the place where the spin condition on $F$ comes in, which guarantees that  
the function $H(t, \tau)$ is defined as the equivariant index of 
an elliptic operator, which is a virtual character of the $S^1$ representation and
therefore is holomorphic for $(t, \tau)\in {\bf R}\times {\bf H}$. 

\vspace{.15in}

{\em Proof of Lemma 4.3.}  Let $z= e^{2 \pi i t}$,  $K= {\rm max}_{j, N} |m_j|$. 
Denote by $D_K\subset \bC ^2$ the domain 
\begin{eqnarray}\begin{array}{l}
|q|^{1/K} < |z| < |q|^{-1/K},\ \ \ \  0< |q| < 1.
\end{array}\nonumber\end{eqnarray}
Let $f_N$ be the contribution of the fixed component $N$ in the function
 $H$. Then by (4.1), (4.4), in $D_K$, $f_N$ has expansion of the form
$$\prod_j \left(1 - z^{m_j}\right)^{-p(p+r)} \sum_{n=0}^\infty
b_{N, n}(z) q^n,$$
where $\Sigma_{n=0}^\infty
b_{N, n}(z) q^n$ is a holomorphic function  of $(z,q)\in D_K$, 
and $b_{N,n}(z)$ are polynomial functions of $z$. So as a meromorphic 
function, in $ D_K$, $H$ has an expansion of the form 
$$\sum_{n=0}^\infty      b_{ n}(z) q^n$$
with each $b_n(z)$ a rational function of $z$, which can only have poles 
on the unit circle $\{ z:|z| = 1\}$.

Now if we multiply the function   $H$  by a function of the form
$$f(z)= \prod_{N} \prod_j \left(1 - z^{m_j}\right)^{p(p+r)},$$
where $N$ runs over the connected components of the fixed point set, 
we get a holomorphic function which has convergent power series expansion
 of the form
$$ \sum_{n=0}^\infty c_n (z) q^n$$
with $\{ c_n (z)\}$  polynomial functions of $z$ in $D_K$. 

By comparing the above two expansions, we get
\begin{eqnarray}\begin{array}{l}
c_n(z) = f(z) b_n(z).
\end{array}\nonumber\end{eqnarray}

On the other hand, we can expand the element 
$\Psi_q(F)\otimes \varphi(F^\bot)$ into 
formal power series of the form $\Sigma_{n=0}^\infty R_n q^n$ with 
$R_n \in K(M)$. 
So, for $t\in [0, 1]\setminus 
 \bQ, z= e ^{2 \pi i t}$, by applying the equivariant index formula 
to each term we get a formal power series of $q$ for 
$H$:
$$\sum_{n=0}^\infty \left(\sum_{m= -N(n)}^{N(n)} a_{m,n} z^m \right) q^n$$
with $a_{m,n}\in {\bC}$ 
and  $N(n)$ some positive integer  depending on $n$.

By comparing the above formulas we get for $t\in [0, 1]\setminus 
 \bQ,\ z= e ^{2 \pi i t}$,
$$b_n (z) = \sum_{m= -N(n)}^{N(n)} a_{m,n} z^m.$$
Since both sides are analytic functions of $z$, this equality holds 
for any $z\in \bC$.

By using the Weierstrass preparation theorem, we then deduce that
$$  \sum_{n=0}^\infty
b_{ n}(z) q^n = {1 \over f(z)}\sum_{n=0}^\infty
 c_n (z) q^n $$
is holomorphic on $(z,q)$ in $D_K$ which clearly contains the set
$\{ (t,q):t\in {\bf R}, q\in {\bf H}\}$.
Q.E.D.\\

We recall that a (meromorphic) Jacobi form  of index $m$ and weight $k$ over 
$L\rtimes \Gamma$, where $L$ is an integral lattice  in the complex plane 
$\bC$ preserved by the modular subgroup $\Gamma \subset SL(2, \bf{Z})$,
 is a (meromorphic) function $F(t, \tau)$ on $\bf{C} \times \bf{H}$ such that
$$F\left({t \over c \tau + d}, { a \tau + b \over c \tau +d}\right) 
= (c \tau + d)^k e^{2 \pi i m ( c t^2 / (c \tau + d))} F(t, \tau),$$
$$F\left(t + \lambda \tau + \mu, \tau\right) = 
e ^{- 2 \pi i m ( \lambda ^2 \tau + 2 \lambda t)} F(t, \tau),\eqno(4.13)$$
where $(\lambda, \mu) \in L$ and $\gamma = \left (\begin{array}{l} a \quad b\\
c \quad d
\end{array}  \right ) \in \Gamma$. If $F$ is holomorphic on $\bf{C} \times \bf{H}$, 
we say that $F$ is a holomorphic Jacobi form [EZ].

\vspace{.15in}

{\bf Lemma 4.5.} {\em The function $H(t, \tau)$ is a holomorphic Jacobi form
of weight $p+r-l$ and index $-n/2$ over $(2\bZ)^2 \rtimes SL(2,\bZ)$.}
 
{\em Proof.} For $\gamma = \left ( \begin{array}{l} a \quad b\\
c \quad d \end{array} \right ) \in SL (2,\bZ)$,  we define its modular 
transformation on $\bC \times \bH$ by
\begin{eqnarray}\begin{array}{l}
\displaystyle{
\gamma (t, \tau) = \left ( { t \over c \tau + d}, {a \tau + b \over c \tau + d}
\right ). }
\end{array}\nonumber\end{eqnarray}
Recall the two generators of $SL(2,\bZ)$ are 
$ S = \left (  \begin{array}{l} 0 \quad -1 \\
1 \qquad    0
 \end{array} \right )\  {\rm and}\ 
T = \left (  \begin{array}{l} 1 \quad 1 \\
0 \quad 1
 \end{array} \right ),$
which act on $\bC \times \bH$ in the following way:
\begin{eqnarray}
S(t,\tau) = \Big ({t \over \tau}, - {1 \over \tau}\Big ),\  \ \ \ 
T(t, \tau) = (t, \tau +1).\nonumber
\end{eqnarray}

So Lemmas 4.1-4.3 imply that $H(t, \tau)$ is a 
(meromorphic) Jacobi form of weight $p+r-l$ and index $-n/2$ over 
$(2\bZ)^2 \rtimes SL(2,\bZ)$.  We now show that it actually is holomorphic. 
 
{}From (4.1) and (4.4)
we know that the possible poles of $H$ in ${\bf C}\times {\bf H}$ are 
of the form 
$$t=\frac{h}{s}(c\tau+d),\eqno(4.14)$$ 
where $h, \ s,\ c, \ d$ are integers with $(c, d)=1$ 
 or $c=1,\ d=0$.

We can always find  integers $a,\ b$ such that $ad-bc=1$. Then 
the matrix 
$$\gamma=\left ( \begin{array}{l}d\quad -b\\
-c \quad a
\end{array} \right ) \in SL(2, {\bf Z})\eqno(4.15)$$ 
induces an action
$$H(\gamma(t, \tau))= H\left(\frac{t}{-c\tau+a}, \frac{d\tau-b}{-c\tau+a}
\right).\eqno(4.16)$$

Now, if $t={h \over s} (c \tau + d)$ is a polar divisor of $H(t,\tau)$, 
then one polar divisor of $H(\gamma(t,\tau))$ is given by
$$\frac{t}{ -c \tau + a}= \frac{h}{s}
\left( c \frac{d \tau -b}{ -c \tau + a} + d \right),\eqno(4.17)
$$
which exactly gives $t = h/s$. But by Lemma 4.1,  up to 
a factor that is holomorphic in $(t, \tau)\in {\bf C}\times {\bf H}$, 
 $H(\gamma(t, \tau))$ is still equal to $H(t, \tau)$ which is holomorphic 
for $t\in {\bf R}$. This implies that $H(t, \tau)$ has no poles in ${\bf C}\times {\bf H}$. 
Q.E.D.

$\ $

{\em Proof of Theorem 3.1}. Since by (4.7), $n=\sum_j (\dim_\sC E_j) m_j^2$, we have $n\geq0$. 

(i) If $n=0$, then since the $S^1$-action is nontrivial, 
it has no fixed point on $M$. Thus
all the Lefschetz number $L(g)$ vanishes by the fixed point formula.

(ii) 
If $n>0$, then Lemma 4.5
shows that $H(t, \tau)$ is a holomorphic Jacobi form of negative index. 
By [EZ, Theorem 1.2], $H(t, \tau)$  must be zero. 

By (i), (ii) and our choice of $\varphi(F^\bot)$, one gets (3.5)
easily. Formula (3.6) then follows from (3.5) and the multiplicativity of the
Witten elements:
$$\Psi_q(TM)=\Psi_q(F)\cdot \Psi_q\left(F^\bot\right).\eqno(4.18)$$

Now for a general $\varphi(F^\bot)$, we write
$$L\left(F^\bot\right){\rm ch}\left(\varphi\left(F^\bot\right)\right)
=\sum_{i}\omega_i\left(F^\bot\right) \ \ {\rm with}\ \ \omega_i\left(F^\bot\right) 
\in H^i(M;{\bf Q}).\eqno(4.19)$$
Then for each $\omega_i(F^\bot) $, one can find an integral polynomial
$\varphi_i(F^\bot)$ such that
$$L\left(F^\bot\right){\rm ch}\left(\varphi_i\left(F^\bot\right)\right)
=n_i\cdot\omega_i\left(F^\bot\right)\eqno(4.20)$$
for some nonzero integer $n_i$. One then verifies that the equivariant
index of $D_{F,\Psi_q(F)\otimes \varphi(F^\bot)}$ can be expressed as a linear
combination with rational coefficients of the equivariant indices of
$D_{F,\Psi_q(F)\otimes \varphi_i(F^\bot)}$'s, which  have been proved to be zero.

The proof of Theorem 3.1 is completed.  Q.E.D.

$\ $

{\em Proof of Corollary 3.2.} We can either use the simple exact sequence 
for the $S^3$-equivariant cohomology groups, 
$$H^4(BS^3)\rightarrow H^4_{S^3}(M) \rightarrow H^4(M)\eqno(4.21)$$
which follows from the spectral sequence for the fibration 
$$ES^3\times_{S^3} M\rightarrow BS^3$$
by using
the fact that $BS^3$ is $3$-connected. Here $ES^3$ is the universal 
$S^3$-principal bundle over the classifying space $BS^3$ of $S^3$. 

Alternatively, one may  prove this by using
 the following simple observation. In fact, at least formally,
 we may write $p_1(F)_{S^1}$ as 
$$p_1(F)_{S^1}= p_1(F) +Au +n\cdot \pi^*u^2\eqno(4.22)$$ 
with $A$ a two form on $M$.

If the $S^1$-action is induced from an $S^3$-action, then $p_1(F)_{S^1}$ 
should be invariant under the Weyl group action 
$u\rightarrow -u$, which implies $A=0$. 

This means under the condition of Corollary 3.2, there exists  $n \in \bZ$,
 such that $p_1(F)_{S^1}= n \pi^* u^2$. By Theorem 3.1, we get Corollary 3.2.
 Q.E.D.

$\ $

By using a special case of the above arguments, one gets the following result
in which the condition on $p_1(F)$ is no longer needed (Compare with [HL1, Prop. 3.2]).
It generalizes the classical Atiyah-Hirzebruch vanishing theorem [AH] to the foliated 
manifolds.

$\ $

{\bf Theorem 4.6.} {\it If $S^1$-acts effectively
on a compact connected foliation $(M,F)$ and preserves
the leaves induced by $F$ as well as the spin structure on $F$, then 
$\langle\widehat{A}(TM),[M]\rangle =0$.}

\vspace{.20in}

{\bf \S 5. Manifolds with split tangent bundle }

\vspace{.15in}

In this section, we no longer assume that $F$ is integrable.
We assume instead that $M$ itself is a spin manifold and the $S^1$-action
preserves the spin structures on $TM$ and $F$. Then it also preserves the
induced spin structure on $TM/F\simeq F^\bot$.  
Consequently, the $S^1$-action on the restriction
of $F^\bot$ to the fixed point set of the $S^1$-action on $M$
need not be trivial.
 
Let us introduce elements
$$R_q\left(F^\bot\right)=\otimes_{j=1}^\infty {\rm Sym}_{q^j}\left(F^\bot-\dim F^\bot\right)
\otimes \left(\otimes_{m=1}^\infty\Lambda_{q^{m}} \left(F^\bot-\dim F^\bot\right)
\right),$$
$$R'_q\left(F^\bot\right)=\otimes_{j=1}^\infty {\rm Sym}_{q^j}\left(F^\bot-\dim F^\bot\right)
\otimes \left(\otimes_{m=1}^\infty\Lambda_{q^{m-1/2}} \left(F^\bot-\dim
 F^\bot\right)\right),$$
$$R''_q\left(F^\bot\right)=\otimes_{j=1}^\infty {\rm Sym}_{q^j}\left(F^\bot-\dim F^\bot\right)
\otimes \left(\otimes_{m=1}^\infty\Lambda_{-q^{m-1/2}}\left(F^\bot-\dim
 F^\bot\right) \right).\eqno(5.1)$$

Recall that since the $S^1$-action preserves $g^{TM}$, it 
preserves the orthogonal splitting 
$$TM =F\oplus F^\bot.\eqno(5.2)$$

Let $D$ denote the canonical Dirac operator on $M$ associated to $g^{TM}$.
We also consider the twisted Dirac operators 
$$D_{\Psi_q(F)\otimes R'_q\left(F^\bot\right)}=D\otimes \Psi_q(F)\otimes 
R'_q\left(F^\bot\right)$$
and
$$D_{\Psi_q(F)\otimes R''_q\left(F^\bot\right)}=D\otimes \Psi_q(F)\otimes 
R''_q\left(F^\bot\right).$$

Under the above assumptions and notations, 
the main result of this section can be stated as follows.

\vspace{.15in}

{\bf Theorem 5.1.} {\em 
If $p_1(F)_{S^1} =n\cdot \pi^* u^2$ for some integer $n\neq 0$, 
then the equivariant 
indices of $D_{F,\Psi_q(F)\otimes R_q(F^\bot), +}$,
$D_{\Psi_q(F)\otimes R'_q(F^\bot), +}$ and
 $D_{\Psi_q(F)\otimes R''_q(F^\bot), +}$  vanish. In particular, the following
three formulas hold,}
$$\left\langle\widehat{A}(F)L\left(F^\bot\right)
{\mbox{\rm ch}}\left(\Psi_q(F)\right){\mbox{\rm ch}}
\left(R_q\left(F^\bot\right)\right), [M]\right\rangle=0,$$
$$\left\langle\widehat{A}(TM)
{\mbox{\rm ch}}\left(\Psi_q(F)\right){\mbox{\rm ch}}
\left(R'_q\left(F^\bot\right)\right), [M]\right\rangle=0,$$
$$\left\langle\widehat{A}(TM)
{\mbox{\rm ch}}\left(\Psi_q(F)\right){\mbox{\rm ch}}
\left(R''_q\left(F^\bot\right)\right), [M]\right\rangle=0.\eqno(5.3)$$

\vspace{.15in}

{\em Proof of Theorem 5.1}. Let  
$$F|_N = F_0 \oplus \left(\oplus_j E_j\right),\ \ 
 \ \left. F^\bot\right|_N= F^\bot_0\oplus \left(\oplus_j L_j\right)\eqno(5.4)$$
be the corresponding equivariant decomopositions of $F$ and $F^\bot$,
when restricted to the 
connected component $N$ of the fixed point set of the $S^1$-action on $M$. 
Assume the generator $g=e^{2\pi it}\in S^1$ acts on $E_j$ by 
multiplication by $e^{2\pi im_jt}$ and on $L_j$ by multiplication by $e^{2\pi i n_jt}$.

 Let $\{2 \pi i x_j^k \}$ denote the  Chern roots  of $E_j$ 
and $\{ 2 \pi i z_j^k \}$ denote the  Chern roots of $L_j$. We also denote by $\{\pm 2 \pi iy_j\}$ 
and $\{ \pm 2 \pi i w_j\}$ the Chern roots of $F_0 \otimes_\sR \bC$ and $F_0^\bot \otimes_\sR \bC$ respectively.

Let $\theta_1(v, \tau)$,  $\theta_2(v, \tau)$ and  $\theta_3(v, \tau)$ be 
the three theta functions (cf. [Ch]):
$$\theta_3(v, \tau)=c(q)\prod_{n=1}^\infty \left(1 + q^{n-1/2} e^{2 \pi i v}\right) 
\prod_{n=1}^\infty \left(1 + q^{n-1/2} e ^{-2 \pi i v}\right),$$
$$\theta_2(v, \tau)=c(q)\prod_{n=1}^\infty \left(1 - q^{n-1/2} e^{2 \pi i v}\right) 
\prod_{n=1}^\infty \left(1 - q^{n-1/2} e ^{-2 \pi i v}\right),$$
$$\theta_1(v, \tau)=c(q) q^{1/8} 2 \cos(\pi  v)
\prod_{n=1}^\infty \left(1 + q^{n} e^{2 \pi i v}\right) 
\prod_{n=1}^\infty \left(1 + q^{n} e ^{-2 \pi i v}\right).$$

Let us write
$$G_0(\tau)= \widehat{A}(F_0)L\left(F_0^\bot\right)
{\mbox{ch}}\left(\Psi_q\left(F_0\right)\right){\mbox{ch}}\left(R_q
\left(F_0^\bot\right)\right),$$
$$G'_0(\tau)= \widehat{A}\left(TN\right){\mbox{ch}}
\left(\Psi_q\left(F_0\right)\right)
{\mbox{\rm ch}}\left(R'_q\left(F_0^\bot\right)\right),$$
$$G''_0(\tau)= \widehat{A}\left(TN\right)
{\mbox{ch}}\left(\Psi_q\left(F_0\right)\right){\mbox{ch}}
\left(R''_q\left(F_0^\bot\right)\right).\eqno(5.5)$$

By applying the equivariant
index formula (2.21), we get three functions, 
$$G(t, \tau) =\sum_N (2 \pi i)^{-(p+r-{\dim N \over 2})} \left\langle
G_0( \tau)\prod_{j,k} \frac{\theta'(0, \tau)}{\theta(x^k_j+m_jt, \tau)}
\prod_{j,k}\frac{\theta_1(z_j^k+n_jt, \tau)\theta'(0, \tau)}
{\theta(z_j^k+n_jt, \tau)
\theta_1(0, \tau)}, [N]\right\rangle,$$
$$G'(t, \tau) =\sum_N  (2 \pi i)^{-(p+r-{\dim N \over 2})} \left\langle
G'_0( \tau)\prod_{j,k} \frac{\theta'(0, \tau)}{\theta(x^k_j+m_jt, \tau)}
\prod_{j,k}\frac{\theta_2(z_j^k+n_jt, \tau)\theta'(0, \tau)}
{\theta(z_j^k+n_jt, 
\tau)\theta_2(0, \tau)}, [N]\right\rangle,$$
$$G''(t, \tau)=\sum_N (2 \pi i)^{-(p+r-{\dim N \over 2})} \left\langle
G''_0( \tau)\prod_{j,k} \frac{\theta'(0, \tau)}{\theta(x^k_j+m_jt, \tau)}
\prod_{j,k}\frac{\theta_3(z_j^k+n_jt, \tau)\theta'(0, \tau)}
{\theta(z_j^k +n_jt, \tau)
\theta_3(0, \tau)}, [N]\right\rangle\eqno(5.6)$$
corresponding to the equivariant indices of the 
three elliptic operators $D_{F, \Psi_q(F)\otimes R(F^\bot), +}$, 
$D_{\Psi_q(F)\otimes  R'(F^\bot),+}$ and 
$D_{ \Psi_q(F)\otimes R''(F^\bot),+}$ respectively.

Now recall the definitions of the following three modular subgroups:
$$\qquad \begin{array}{l}
\Gamma_0(2) = \left \{ \left ( \begin{array}{l} a \quad b\\
c \quad d  \end{array}\right ) \in SL(2, {\bf Z}): c \equiv 0\ (\rm mod\ 2)
 \right \}, \\
\Gamma^0(2) = \left \{ \left ( \begin{array}{l} a \quad b\\
c \quad d \end{array}\right ) \in SL( 2, {\bf Z}): b \equiv 0\ (\rm mod\ 2) 
\right \},\\
\Gamma_{\theta} = \left \{ \left ( \begin{array}{l} a \quad b\\
c \quad d \end{array} \right ) \in SL( 2, \bf{Z}): 
\left ( \begin{array}{l} a \quad b\\
c \quad d \end{array} \right )
 \equiv   \left ( \begin{array}{l} 1 \quad 0\\
0 \quad 1 \end{array} \right )\ {\rm  or }\ 
\left ( \begin{array}{l} 0 \quad 1 \\
1 \quad 0 \end{array} \right )\ (\rm mod\ 2) 
\right \}.
\end{array}$$

By using the modular transformation formula of the theta-functions [Ch], we 
can immediately prove the following result by proceeding as in the proofs
of Lemmas 4.1 and 4.2.

\vspace{.15in}

{\bf Lemma 5.2.} {\em If $p_1(F)_{S^1}= n\cdot \pi^* u^2$, then $G(t, \tau)$ is 
a Jacobi form over $(2{\bf Z})\rtimes \Gamma_0(2)$, $G'(t, \tau)$ is 
a Jacobi form over  $(2{\bf Z})\rtimes \Gamma^0(2)$ and  $G''(t, \tau)$ is 
a Jacobi form over  $(2{\bf Z})\rtimes \Gamma_\theta$.}
All of them are of index $-n/2$ and weight $p+r$.
\vspace{.15in}

For $\gamma=  \left ( \begin{array}{l} a \quad b\\
c \quad d \end{array}\right ) \in SL(2, \bf{Z})$, let us use the notation
$$H(\gamma(t,\tau))|_{m,k} = (c \tau + d)^{-k} e^{- 2 \pi i m c t^2/(c \tau +d)} 
H\left({t \over c \tau + d}, { a \tau +b \over  c\tau +d }\right)\eqno(5.7)$$
to denote the action of $\gamma$ on a Jacobi form $H$ 
of index $m$ and weight $k$.

Recall that  $2p=\dim F$ and $2r=\dim F^\bot$.
The following lemma can be proved easily by proceeding as the proof of 
Lemma 4.3.  

\vspace{.15in}

{\bf Lemma 5.3.}  {\em  For any $\gamma\in SL(2, \bf{Z})$, let $F(t, \tau)$ be one of the 
functions 
$G(t,\tau)$, $G'(t,\tau)$ and $G''(t,\tau)$. Then $F(\gamma(t, \tau))|_{{-n \over 2}, p+r}$
is holomorphic for $(t, \tau)\in \bf{R} \times\bf{H}$.}

\vspace{.15in}

Again this is the place where the index theory comes in to cancel part 
of the poles of these  functions. Here the spin  conditions on $F$, $TM$
are crucially needed.

Now, by using the same argument as in the proof of Lemma 4.5,  we get

\vspace{.15in}

{\bf Lemma 5.4.} {\em For a (meromorphic) Jacobi form $H(t, \tau)$ of index $m$ and 
weight $k$ over $L\rtimes \Gamma$, assume that $H$ may only have polar 
divisors of the form $t=(c \tau + d)/l$ in $\bf{C} \times \bf{H}$ for some 
integers $c,d$ and $l\neq 0$. If $H(\gamma(t, \tau))|_{m,k}$ is holomorphic 
for $t\in \bf{R}$, $\tau \in \bf{H}$ for every $\gamma\in SL(2, \bf{Z})$, 
then $H(t, \tau)$ is holomorphic for any $t\in \bf{C}$ and $\tau \in \bf{H}$.}

\vspace{.15in}

{}From Lemmas 5.3 and 5.4 one sees, as in the proof of Theorem 3.1, that the 
$G$'s are holomorphic Jacobi forms of index $-n/2$, and therefore must be 
zero. (Here we have used the fact that
$n>0$.)

The proof of Theorem 5.1 is completed.  Q.E.D.

\vspace{.15in}

{\bf Remark 5.5.} If $n=0$, then we get the rigidity properties in Theorem 5.1 
instead of the vanishing results.

$\ $

{\bf \S 6. Concluding remarks.}

\vspace{.15in}

Motivitate by  Corollary 3.2, we find it is 
interesting and reasonable to make the following conjecture
which may be viewed as a foliation analogue of a conjecture of Hoehn and 
Stolz [S].

We consider an oriented compact foliation $M$ which is foliated by 
a spin integrable subbundle $F$ of $TM$. 
Let $g^F$ be a metric on $F$.
 
\vspace{.15in}

{\bf Conjecture 6.1.} {\em If ${1 \over 2} p_1(F)=0$, and the Ricci 
curvature of $g^F$ along 
each leaf is positive, then the Witten genus of $M$, 
$\langle\widehat{A}(TM){\mbox{\rm ch}}(\Psi_q(TM)), [M]\rangle$, vanishes. }

\vspace{.15in}

As have been remarked in the introduction, we may as well 
take $\psi(F)$ or $\varphi(F^\bot)$ in Section 2 as the elements in 
$K(M)$ induced from loop group representations. Then 
the modularity of the characters of the loop group representations can be used
to prove vanishing theorems for the correponding
twisted sub-Dirac operators. 
On the other hand, the construction of the sub-elliptic operators is very flexible. 
For example if the the integrable subbundle of the foliation has almost complex structure or 
$\mbox{Spin}^c$-structure, then we can construct sub $\bar{\partial}$-operator 
or $\mbox{Spin}^c$ sub-Dirac operator correspondingly. 
If there exists a compact Lie group action on $M$ preserving the 
leaves, the the rigidity and vanishing theorems can be proved for the equivariant indices 
of such operators which generalize the corresponding rigidity and vanishing 
results for the usual elliptic genera. See [Liu2] or [LiuMa2] for some details about these. 

In concluding, we may also replace the signature operator in the normal direction by 
other elliptic operators like the de Rham type operator from which we can 
derive the vanishing of characteristic numbers like 
$$\left\langle\widehat{A}(TM){\mbox{ ch}}\left(\Psi_q(F)
\right) e\left(F^\bot\right), 
[M]\right\rangle,$$
where $e(F^\bot)$ denotes the Euler class of $F^\bot$.

 {\small   

$\ $

Kefeng Liu,
Department of Mathematics, Stanford University, Stanford, CA 94305, USA.

{\em E-mail address}: kefeng@math.stanford.edu

$\ $

Xiaonan Ma,
Humboldt-Universitat zu Berlin, Institut f\"ur Mathematik, unter den Linden 6,
D-10099 Berlin, Germany.

{\em E-mail address}: xiaonan@mathematik.hu-berlin.de

$\ $

Weiping Zhang,
Nankai Institute of Mathematics, 
Tianjin 300071, P. R. China.}

{\em E-mail address}: weiping@nankai.edu.cn


\begin {thebibliography}{15}

\bibitem [AH]{}  Atiyah M.F. and  Hirzebruch F., Spin manifolds and groups 
actions, {\it Essays on topology and Related Topics, Memoires d\'edi\'e
\`a Georges de Rham} (ed. A. Haefliger and R. Narasimhan),
Springer-Verlag, New York-Berlin (1970), 18-28.


\bibitem [AS]{} Atiyah M.F., Singer I.M., The index of elliptic operators III. 
{\em Ann. of Math}. 87 (1968), 546-604.

\bibitem [BT]{} Bott R. and  Taubes C., On the rigidity theorems of Witten, 
{\em J.A.M.S}. 2 (1989), 137-186.

\bibitem [Ch]{} Chandrasekharan K., {\em Elliptic functions}, Springer,
 Berlin (1985).

\bibitem [Co]{} Connes A., {\em Noncommutative Geometry}, Academic Press, 1994.

\bibitem [EZ]{} Eichler M., and Zagier D., {\em The theory of Jacobi forms},
 Birkhauser, Basel, 1985.

\bibitem [HL1]{} Heitsch, J., Lazarov C., A Lefschetz fixed point theorem 
for foliated manifolds, {\em Topology} 29 (1990), 127-162.

\bibitem [HL2]{} Heitsch, J., Lazarov C., Rigidity theorems for foliations by 
surfaces and spin manifolds, {\em Michigan Math. J.} 38 (1991), 285--297.


\bibitem [Liu1]{}  Liu K., On elliptic genera and theta-functions, 
{\em Topology} 35 (1996), 617-640.

\bibitem [Liu2]{}  Liu K., On modular invariance and rigidity theorems, 
{\em J. Diff. Geom}. 41 (1995), 343-396.

\bibitem [LiuMa1]{}  Liu K., Ma, X., On family rigidity theorems I, to appear 
in {\em Duke Math. J.}.

\bibitem [LiuMa2]{}  Liu K., Ma, X., On family rigidity theorems II, 
{\em preprint}.

\bibitem [LiuZ]{}  Liu K., Zhang W., Adiabatic limits and foliations, {\em preprint}.  

\bibitem [S]{} Stolz, S., A conjecture concerning positive Ricci curvature and the
Witten genus. {\em Math. Ann.} 304 (1996), 785-800.

\bibitem [T]{} Taubes C., $S^1$-actions and elliptic genera, 
 {\em  Comm. Math. Phys.} 122 (1989), 455-526.

\bibitem [W]{} Witten E., The index of the Dirac operator in loop space, 
in {\em Elliptic Curves and Modular forms in Algebraic Topology},
 Landweber P.S. eds., SLNM 1326, Springer, Berlin, 161-186.
\end{thebibliography}
\end{document}